\documentclass[11pt]{amsart}
\usepackage{amssymb}
\numberwithin{equation}{section}
     \newtheorem{thm}{Theorem}[section]

     \newtheorem{lem}[thm]{Lemma}
\theoremstyle{definition}
     \newtheorem{defn}{Definition}[section]
     
\theoremstyle{remark}
     \newtheorem{rem}{Remark}[section]

\def\l{\lambda}

\def\s{\sigma}

\def\d{\delta}
\def\D{\Delta}

\def\mc{\mathcal}

\def\G{\Gamma}
\def\e{\epsilon}

\def\O{\Omega}

\def\al{\alpha}
\def\bt{\beta}
\def\ro{\rho}

\def\ra{\rightarrow}
\def\nd{\noindent}
\def\p{\partial}

\def\la{\langle}
\def\ra{\rangle}

\def\vp{\varphi}

\begin{document}

\author[Bin Xu]
{Bin Xu}
\title[asymptotic behavior of eigenfunctions]
{Asymptotic behavior of $L^2$-normalized eigenfunctions of the Laplace-Beltrami operator
on a closed Riemannian manifold}

\address{Bin Xu\\
Nankai Institute of Mathematics\\
Nankai University\\
Tianjin 300071 China (The former addresses: Department of Mathematics\\
Tokyo Institute of Technology\\
2-12-1 Oh-okayama, Meguro-ku\\
Tokyo 152-8551 Japan; Department of Mathematics\\ Johns Hopkins University\\
 Baltimore, MD 21218 USA)
}
\email{bxu@math.jhu.edu}

\subjclass[2000]{35P20, 35L05.}
\keywords{Laplace-Beltrami operator, unit band spectral projection operator, eigenfunctions.}

\begin{abstract}
Let $e(x,y,\l)$ be the spectral function and ${\chi}_\l$ the unit band
spectral projection operator,  with respect to the Laplace-Beltrami operator $\D_M$
on a closed Riemannian manifold $M$. We firstly review
the one-term asymptotic formula of $e(x,x,\l)$ as $\l\to\infty$ by H{\" o}rmander (1968)
 and the one of $\p^\al_x\p^\bt_y e(x,y,\l)|_{x=y}$ as $\l\to\infty$ in a
geodesic normal coordinate chart by the author (2004)
and the sharp asymptotic estimates from above of the mapping norm $\|\chi_\l\|_{L_2\to L_p}$
($2\leq p\leq\infty$) by Sogge (1988 $\&$ 1989) and of the
mapping norm  $\|\chi_\l\|_{L_2\to {\rm Sobolev}\ L_p}$
by the author (2004).
In the paper we show the one term asymptotic formula for
$e(x,y,\l)$ as $\l\to\infty$, provided that the Riemannian distance between $x$ and $y$ is ${\rm O}(1/\l)$.
As a consequence, we obtain the sharp estimate of the mapping norm
$\|\chi_\l\|_{L_2\to C^\d}$ ($0<\d<1$), where $C^\d(M)$ is the space of H{\" o}lder continuous functions
with exponent $\d$ on $M$.
Moreover, we show a geometric property of the eigenfunction $e_\l$: $\D_M e_\l+\l^2 e_\l=0$, which says that
$1/\l$ is comparable to the distance between the nodal set of $e_\l$ (where $e_\l$ vanishes)
and the concentrating set of $e_\l$ (where $e_\l$ attains its maximum or minimum) as $\l\to\infty$.

\end{abstract}

\maketitle
\tableofcontents
%%%==========================================================
\section{Introduction}
R. Balian and C. Bloch \cite{BB1} \cite{BB2} \cite{BB3} studied thoroughly the distribution of eigenvalues
of the wave equation
\begin{equation}
\label{equ:Hel}
\D\varphi+k^2 \varphi=0
\end{equation}
in a volume $V$ with smooth surface $S$, where $\D$ is the Euclidean Laplacian and $\varphi$ satisfies
some boundary condition. For instance, the wave function (eigenfunction) $\varphi$ may represent a sound wave, either in
dimension three (room acoustics), or in dimension two \cite{P1} \cite{P2} (vibration of a membrane),
and then the eigenvalue $k^2$ represents the frequency of $\varphi$. On the other hand,
at the end of the eighteenth century Ernst Florens Chladni conceived an experiment which can be considered as a
precursor of the research of waves. \\

\nd
{\small
Chladni had an inspired idea:  to make sound waves in a solid material visible.
This he did by getting metal plates to vibrate, stroking them with a violin bow.
Sand or a similar substance spread on the surface of the plate naturally
settles to the places where the metal vibrates the least,
making such places visible. These places are the so-called nodes, which are wavy lines on the surface.
The plates vibrate at pure, audible pitches, and each pitch has a unique nodal pattern.
Chladni took the trouble to carefully diagram the patterns, which helped to popularize his work.
Then he hit the lecture circuit, fascinating audiences in Europe with live demonstrations. This culminated with a
command performance for Napoleon, who was so impressed that he offered a prize to anyone who could explain the patterns.
More than that, according to Chladni himself, Napoleon remarked that irregularly shaped plate would be much harder to
understand! While this was surely also known to Chladni, it is remarkable that Napoleon had this insight. Chladni
received a sum of 6000 francs from Napoleon, who also offered 3000 francs to anyone who could explain the patterns.
The mathematician Sophie Germain took the prize in 1816, although her solutions were not completed until the work of
Kirchoff thirty years later. Even so, the patterns for irregular shapes remained (and to some extent remains) unexplained.
(The author took the paragraph from the website of Eric. J. Heller Gallery (http://www.ericjhellgallery.com).
The readers can also enjoy Chladni's beautiful diagram by accessing {\it Quantum modes and Classical Analog} in the
gallery.)
}

Actually the vibrations of stiff plates in these 200 year old experiments can be described by the equation
(Cf $\S\,$25 of \cite{Lan}):
\begin{equation}
\label{equ:Lan}
(\D)^2\psi_j=k_j^4 \psi_j,
\end{equation}
where $\D=\p^2/\p x^2+\p^2/\p y^2$ is the two dimensional Laplacian, $\psi_j(x,\,y)$ is the amplitude function
of the $j$th resonance, and $k_j$ the associated frequency. Note that in contrast to the vibrations of membranes without
internal stiffness it is the square of the Laplacian which enters \eqref{equ:Lan}. Thus the functions obeying the
Helmholtz equation \eqref{equ:Hel} form a class of solutions of \eqref{equ:Lan}. If the plates are clamped along
the outer rim, the problem is completely equivalent to the quantum mechanics of a particle in a box with
infinitely high walls (the state of the particle satisfies the Dirichlet boundary condition).

A mathematical version of the above problem in physics is as follows:
on a compact Riemannian manifold $M$ without boundary, describe the asymptotic behaviour of the eigenfunctions of the
Laplace-Beltrami operator $\D_M$ on $M$ in terms of the global geometry of $M$, as the associated eigenvalues tend to infinity.
However, this mathematical problem remains yet unsolved but a lot of
partial results have been obtained by Sarnak \cite{Sar}, Sogge \cite{S1} \cite{S2},  Sogge-Zelditch \cite{SZ},
Toth-Zelditch \cite{TZ} and Zelditch \cite{Z1} et al.
We can only mention the detail of the results in \cite{S1} \cite{S2} because the interesting results
in \cite{Sar} \cite{SZ} \cite{TZ} \cite{Z1} have no direct connection with our story.

The philosophy of the problem we are going to consider is as follows:
thinking of $L_2$-normalized eigenfunctions $e_\l$:
\begin{equation}
\label{equ:eigf}
\D_M e_\l+\l^2 e_\l=0,\ \int_M |e_\l|^2\,{\rm d}x=1,
\end{equation}
as elements in some function spaces $A$ on $M$,
we obtain the upper bounds of $\|e_\l\|_A$ as the associated eigenvalues $\l^2$ tends to $\infty$.
The function spaces $A$ we are going to consider are: $L_p$ ($2\leq p<\infty$) spaces,
fractional Sobolev $L_p$ spaces, $C^k$ space with $k$ non-negative integers and the H{\" o}lder-Zygmund space
(cf \cite{SS} \cite{Tr} for their definitions).
The size of the $L_\infty$ norm, may be considered as a rough measure for how
unevenly the eigenfunction $e_\l$ is distributed over $M$:
If $e_\l$ is "small" (say ${\rm o}(1)$ as $\l\to\infty$) on a "large" set $S$ (say, of area
$|S|=|M|-\d$), then the normalizing condition in \eqref{equ:eigf}  forces $\|e_\l\|_{L_\infty}$ to be at least of order
$1/\sqrt{\d}$; therefore upper bounds on $\|e_\l\|_{L_\infty}$ imply lower bounds on $\d$, i.e. the area of
the set where $u$ is "concentrated". Also, upper bounds on $\|e_\l\|_{L_\infty}$ yield upper bounds
for the multiplicities of the eigenvalue $\l$ (cf Proposition 3 in Grieser \cite{G}).
The sizes of the $C^k$ and the H{\" o}lder norms may be considered as a rough measure for how rapidly
the eigenfunction $e_\l$ oscillates on $M$: By the mean value theorem,
the ratio $\|e_\l\|_{L_\infty}/\|e_\l\|_{C^1}$ bounds from below
 the distance\footnote{As mentioned in the abstract, the distance is actually comparable to $1/\l$ as $\l\to\infty$, which
will be proved in Theorem \ref{thm:nodal}. For safety, we remind
the readers that $A(\l)>0$ is comparable to $B(\l)>0$ means
$C_1B(\l)\leq A(\l)\leq C_2B(\l)$ for two absolute constants
$C_1,\,C_2$. } between its nodal set (the set where $e_\l$
vanishes) and its {\it concentrating set} (where $e_\l$ achieves
its maximum or minimum). The sizes of $L_p$ norms, $2\leq p\leq
\infty$, have deep relationship with the $L_p$ convergence of
Riesz means on $M$ (cf Sogge \cite{S3}). It turns out that these
upper bounds have the expression ${\rm O}(\l^D)$ as $\l\to\infty$,
where $D>0$ are the universal constants depending on $\dim\,M$ and
the parameters of the function spaces $A$, but not depending on
the geometry of $M$. Moreover, we also find some sequences of
spherical harmonics which attain these upper bounds as their
eigenvalues go to infinity.

\section {The $C^k$ norm of $\chi_\l$}

\subsection{The spectral function and the unit band spectral projection operator}
We firstly fix the notations of this paper. Let $(M,\ g)$ be a
smooth closed Riemannian manifold of dimension $n\geq 2$ and
$\Delta$ the positive Laplace-Beltrami operator on $M$. Let
$L_2(M)$ be the space of square integrable functions on $M$ with
respect to the Riemannian density ${\rm d}v(x)$ Let
$e_1(x),\,e_2(x),\,\cdots$ be a complete orthonormal basis in
$L_2(M)$ for the eigenfunctions of $\D$ such that $0\leq
\l_1^2\leq \l_2^2\leq \,\cdots$ for the corresponding eigenvalues,
where $e_j(x)$ ($j=1,2,\cdots$) are real-valued smooth function on
$M$ and $\l_j$ are nonnegative real numbers. Let $\l$ be a
positive real number $\geq  1$ throughout this paper. We define
the {\it spectral function} $e(x,y,\l)$ to be the Schwartz kernel
of the orthogonal projection of $L_2(M)$ onto the subspace spanned
by the eigenfunctions whose eigenvalues lie in $[0,\,\l^2]$. We
define the {\it unit band spectral projection operator} (UBSPO)
$\chi_\l$ to be the orthogonal projection of $L_2(M)$ onto the
subspace spanned by the eigenfunctions whose eigenvalues lie in
$(\l^2,\,(\l+1)^2]$ (cf in \cite{Xu2} a detailed survey for how
Sogge invented the operator $\chi_\l$). It is clear that
$e(x,y,\l)$ is smooth with respect to $(x,\,y)$ in $M\times M$ and
the following relations among $\{e_j\}$, $e(x,y,\l)$ and $\chi_\l$
hold:
\begin{gather}
\label{equ:spect}
e(x,y,\l)=\sum_{\l_j\leq \l} e_j(x) e_j(y), \\
\label{ubs}
{\chi}_\l f(x)=\sum_{\l_j\in (\l,\,\l+1]}\int_M e_j(x) e_j(y)f(y)\,{\rm d}v(y).
\end{gather}

The wave kernel method applied in \cite{A} \cite{G} \cite{H1}
\cite{S2} shows that the spectral function $e(x,y,\l)$ and the
UBSPO $\chi_\l$ turn out to be easier to study than a single
eigenfunction $e_\l$ because they have a local character. We shall
explain the exact meaning of the local character in
Subsection\,\ref{subsec:wk}. To obtain the upper bound of the norm
$\|e_\l\|_A$ for a single eigenfunction $e_\l$ in
\eqref{equ:eigf}, we need only to estimate from above the mapping
norm
$$\|\chi_\l\|_{L_2\to A}=\sup_{\|f\|_{L_2}=1} \|\chi_\l f\|_A$$
which automatically dominates $\|e_\l\|_A$. It is the trick in this paper!

Let $f(\l)$ be a function defined in a neighbourhood of $+\infty$.
Suppose that there exists an estimate $f(\l)={\rm O}(\l^E)$ as $\l\to+\infty$ with
$E$ a constant independent of $\l$. We say the estimate is {\it sharp} if
\begin{equation}
\label{equ:sharp}
\limsup_{\l\to\infty} \frac{|f(\l)|}{\l^E}>0,
\end{equation}
that is, $f(\l)\not={\rm o}(\l^E)$ as $\l\to\infty$.
Following Hardy, we use the notion $f(\l)=\O (\l^E)$ as $\l\to\infty$ to mean that
$f(\l)$ is ${\rm O}(\l^E)$ but not ${\rm o}(\l^E)$ as $\l\to\infty$.
We will consider functions of $\l$ such as $\|\chi_\l\|_{L_2\to A}$,
\begin{equation*}
\sum_{\l_j\in (\l,\,\l+1]} |Pe_j(x)|^2\ \ {\rm and}
\sum_{\l_j\in (\l,\,\l+1]} |Pe_j(x)-Pe_j(y)|^2
\end{equation*}
with $P$ some kind of differential operators on $M$.

For simplicity, we use the notation $\|\cdot\|_p$, $2\leq
p\leq\infty$, to denote the $L_p$ norm on $M$. One can derive the
estimates of $\|\chi_\l\|_{L_2\to L_\infty}$ and of the $L_\infty$
norms of eigenfunctions directly from the wave kernel method (cf
Grieser \cite{G}). However, to obtain more information and to be
faithful to the history of the mathematical story, we shall
introduce the (local) Weyl law, which implies the sharp estimate
of $\|\chi_\l\|_{L_2\to L_\infty}$.

\subsection{The local Weyl law and its generalizations}
Levitan \cite{Lev},
Avakumovi${\rm {\check c}}$ \cite{A} and H{\" o}rmander \cite{H1} obtained a
one-term asymptotic formula
of the spectral function $e(x,x,\l)$ as following:
\begin{equation}
\label{equ:asy}
e(x,x,\l)=\frac{\l^n}{2^n \pi^{n/2}\G(1+\frac{n}{2})} +{\rm O}(\l^{n-1}),\ \l\to\infty.
\end{equation}
Since the integration of \eqref{equ:asy} on $M$ gives the Weyl law,
\[N(\l)=\sharp\{j:\,\l_j\leq \l\}=\frac{|M|}{2^n \pi^{n/2}\G(1+\frac{n}{2})}\l^n+{\rm O}(\l^{n-1}),\]
the one term asymptotic formula \eqref{equ:asy} is called the local Weyl law.
As a consequence, there holds the uniform estimate of a collection of eigenfunctions,
\begin{equation}
\label{equ:spe}
\sum_{\l_j\in (\l,\,\l+1]} |e_j(x)|^2 =\O(\l^{n-1})\ {\rm as}\ \l\to\infty,\  \forall\, x\in M.
\end{equation}
Sogge noted in (2.1-2) of \cite{S4} (implicitly in (1.7) of
\cite{S3}) that \eqref{equ:spe} is equivalent to $\chi_\l$
\begin{equation}
\label{equ:c0}
||\chi_\l f||_{L_2\to L_\infty}=\O (\l^{(n-1)/2})\ {\rm as}\ \l\to\infty.
\end{equation}

Let $k$ be a non-negative integer and
$H_\infty^k(M)$ be the space of functions on $M$ with $L_\infty$ derivatives up to order $k$.
Yuri Safarov (\cite{Sa} Theorem 1.8.5) and the author
(\cite{Xu} Theorem 1.1) generalized the asymptotic formula \eqref{equ:asy} for
the restrictions of the spectral function on the diagonal to those for
the restrictions of derivatives of the spectral function on the diagonal. In particular, the
author proved the following: In a geodesic normal coordinate chart $(X\,,x)$ with sufficiently small
radius\footnote{The radius is actually comparable to the injectivity radius ${\rm inj}_M$ of $M$.} of $M$,
for any two multi-indices $\al,\bt\in {\Bbb Z}_+^n$ the following estimates hold uniformly for $x\in X$
as $\l\to \infty${\rm :}
\begin{equation}
\label{equ:deri}
\p_x^\al\p_y^\bt e(x,y,\l)|_{x=y}=
\left\{
\begin{array}{rl}
C_{n,\al,\bt}\,\l^{n+|\al+\bt|}+{\rm O}(\l^{n+|\al+\bt|-1}) & {\rm if}\
\al\equiv\bt\  {\rm (2)}^\dagger \\
{\rm O}(\l^{n+|\al+\bt|-1}) &   {\rm otherwise},
\end{array}
\right.
\end{equation}
where for multi-indices $\al,\bt$ such that $\al\equiv \bt$ {\rm (2)}\footnote{$^\dagger$
$\al\equiv\bt\ (2)$ means that $\al_j-\bt_j$, $1\leq j\leq n$, are even integers.},
\begin{eqnarray}
\label{equ:const1}
C_{n,\al,\bt}&=&(2\pi)^{-n}\,(-1)^{(|\al|-|\bt|)/2}\,
\int_{B_n}\, x^{\al+\bt}\,dx
\nonumber
\\
&=&
\label{equ:skb}
(-1)^{(|\al|-|\bt|)/2}\frac{\prod_{j=1}^n (\al_j+\bt_j-1)!!}
{\pi^{n/2}\,2^{n+|\al+\bt|/2}\,\G(\frac{|\al+\bt|+n}{2}+1)}
\nonumber
\end{eqnarray}
In particular, if $\al=\bt$, then the following estimate holds uniformly
for $x\in X$ as $\l\to\infty${\rm :}
\begin{equation}
\label{equ:sym}
\sum_{\l_j\leq\l} |\p^\al e_j(x)|^2 =C_{n,\al} \l^{n+2|\al|} +
{\rm O}(\l^{n+2|\al|-1})\ ,
\end{equation}
where $C_{n,\al}=C_{n,\al,\al}>0$.
As a consequence, there hold the following uniform estimates of a collection of derivatives of eigenfunctions,
\begin{equation}
\label{equ:speal}
\sum_{\l_j\in (\l,\,\l+1]} |\p^\al e_j(x)|^2 =\O(\l^{2|\al|+n-1})\ {\rm as}\ \l\to\infty,\  \forall\, x\in M.
\end{equation}
Covering $M$ by a finite number of normal charts, similarly to the equivalence between \eqref{equ:spe} and \eqref{equ:c0},
we can see easily that the validity of the estimates \eqref{equ:speal} for all multi-indices $\al$ with $|\al|\leq k$ is
equivalent to
\begin{equation}
\label{equ:ck}
\|\chi_\l \|_{L_2\to H^k_\infty}=\O (\l^{k+(n-1)/2})\ {\rm as}\ \l\to\infty.
\end{equation}

\begin{rem}
Besides the $\al=\bt$ case in our generalized local Weyl law \eqref{equ:deri}, the other cases also have their own
interest because the argument by Sobolev's embedding theorem only gives a rough estimate
(cf Theorem 17.5.3 in \cite{H2}): in a normal chart, there holds
\begin{equation}
\label{equ:albt}
\p_x^\al \p_y^\bt e(x,y,\l)={\rm O}(\l^{n+|\al+\bt|})\ {\rm as}\ \l\to\infty.
\end{equation}
\end{rem}

\begin{rem}
Grieser \cite{G} and Sogge \cite{S4} showed that \eqref{equ:c0} also holds for both the Dirichlet and the Neumann
Laplacians on a compact Riemannian manifold $N$ with non-empty smooth boundary. Moreover,
X.-J. Xu \cite{X} proved Case $k=1$ of \eqref{equ:ck} on such an $N$ under the same boundary conditions.
\end{rem}

\subsection{The wave kernel method}
\label{subsec:wk}
Let $x$ and $y$ be two points in $M$ and ${\rm dist}\,(x,\,y)$ be the Riemannian distance between $x$ and $y$.
For each fixed $y\in M$, the wave kernel
\begin{equation}
\label{equ:wk}
K(t,x,y)=K_M(t,x,y)=\sum_j \cos\,(t\l_j) e_j(x)e_j(y)
\end{equation}
on $M$ is the solution of
\begin{eqnarray}
\Bigl(\frac{\p^2}{\p t^2}-\D_x\Bigr)K(t,x,y)&=&0\ \ {\rm in}\ {\Bbb R}_t\times M
\nonumber\\
K(0,x,y)&=&\d_y(x)\nonumber\\
\frac{\p}{\p t}K(t,x,y)|_{t=0}&=&0.
\end{eqnarray}
The convergence of \eqref{equ:wk} should be considered in the sense of distribution in $t$, for each fixed
$x$ and $y$. By the finite propagation speed of solutions of the wave equation
(cf H{\" o}rmander \cite{H2}, Lemma 17.5.12,
for example, for an easy proof by energy estimates), we have that
\[{\rm supp}\,K\subset\{(t,x,y):\,{\rm dist}\,(x,y)\leq |t|\}\]
and that $K(t,x,y)$ depends only on the data of $(M,\,g)$ in
\[B_t(x,y):=\{z\in M:\,{\rm dist}\,(x,z)+{\rm dist}\,(y,z)\leq |t|\},
\]
which may thought of as the local character of the wave kernel.
If $|t|\leq {\rm inj}_M$, then we can obtain the approximation of
the wave kernel as precise as we desire by the virtue of the Hadamard parametrix
(cf H{\" o}rmander \cite{H2}, $\S\,$17.4). Moreover, the approximating function of
the wave kernel $K(t,x,y)$ and its derivatives $\p_x^\al \p_y^\bt K(t,x,y)$, as $|t|\leq {\rm inj}_M$,
have explicit expression in normal charts of $M$.
These properties of the wave kernel imply the local character of $e(x,y,\l)$ and $\chi_\l$ mentioned
before. Roughly speaking,
the asymptotic property of $\p_x^\al \p_y^\bt e(x,y,\l)|_{x=y}$ as $\l\to +\infty$
is determined by the inverse Fourier transformation with respect to $t$ of
$\p_x^\al \p_y^\bt K(t,x,y)$ times a cutoff function supported near the origin of ${\Bbb R}_t$
(cf \cite{Xu} for the detail); the UBSPO $\chi_\l$ has the integral kernel
essentially equal to a $C^\infty$ function with support in a neighbourhood of the diagonal in $M\times M$
times
\[\l^{(n-1)/2} \exp\Bigl( \sqrt{-1}\l\,{\rm dist}\,(x,y)\Bigr) \Bigl({\rm dist}\,(x,y)\Bigr)^{(1-n)/2}\]
plus a similar term where $e^{\sqrt{-1}\l\,{\rm dist}\,(x,y)}$ is replaced
by $e^{-\sqrt{-1}\l\,{\rm dist}\,(x,y)}$ (cf Sogge \cite{S2}). We do not want to repeat
the exact statement for the local character property of $\p_x^\al \p_y^\bt e(x,y,\l)|_{x=y}$ and $\chi_\l$ and
the proof, because they need much more pages. However, it is easy for us to explain the local character property
of $\sum_{\l_j\in (\l,\,\l+1]} |e_j(x)|^2$,
which is the square-sum of a bunch of eigenfunctions in general.
The following argument is due to Grieser \cite{G}.

Take $\e={\rm inj}_M$ and choose
a Schwarz function $\ro$ on ${\Bbb R}$ such that
\[\ro>0,\ \ro|_{[0,\,1]}\geq 1,\ {\rm supp}\,{\hat \ro}\subset (-\e,\,\e),\]
where ${\hat \ro}(t)=\int_{-\infty}^\infty e^{-\sqrt{-1}t\l}\ro(\l)\,{\rm d}\l$ denotes the Fourier transform.
See a construction for such a $\ro$ in $\S\,$17.5 \cite{H2}. Consider the sum, convergent by \eqref{equ:albt},
\begin{equation}
\label{equ:ro}
\sum_j \ro(\l-\l_j) e_j(x) e_j(y).
\end{equation}
Writing, by the Fourier inversion formula,
\[\ro(\l-\l_j)=2\bigl({\hat \ro}(t)\cos\,(t\l_j)
\bigr)^\vee-\ro(\l+\l_j)\]
(the superscript $^\vee$ denotes the inverse Fourier transform $t\to \l$), then we obtain
\begin{equation*}
\sum_j \ro(\l-\l_j) e_j(x) e_j(y)=2\bigl(\ro(t)K(t,x,y)\bigr)^\vee(\l)+{\rm O}(\l^{-\infty})
\end{equation*}
where the error term is small by \eqref{equ:ro} and the rapid decay of $\ro$.
This tells us that the weighted sum \eqref{equ:ro} over many eigenfunctions has a local character.
From the assumption on $\ro$ we get
\[\sum_{\l_j\in (\l,\,\l+1]} |e_j(x)|^2\leq 2\bigl(\ro(t)K(t,x,x)\bigr)^\vee(\l)+{\rm O}(\l^{-\infty}).\]
From the exact information on the singularity of $K(t,x,x)$ as $|t|\leq \e$, we can show the estimate,
\[\int_{-\infty}^\infty e^{-\sqrt{-1}t\l} {\hat \ro}(t) K(t,x,x)\,{\rm d}t={\rm O}(\l^{n-1})
\ {\rm as}\ \l\to\infty,\]
which gives an alternative proof for
\[\sum_{\l_j\in (\l,\,\l+1]} |e_j(x)|^2={\rm O}(\l^{n-1})\ {\rm and}\
\|\chi_\l\|_{L_2\to L_\infty}={\rm O}(\l^{(n-1)/2})\]
as $\l\to\infty$.
We remark that the similar weighted sums as \eqref{equ:ro} over derivatives of many eigenfunctions also have
a local character.

\section{The H{\" o}lder-Zygmund norm of $\chi_\l$}
\subsection{Interpolation method} Let $\s$ be a non-negative real number and ${\mc C}^\s (M)$ be
the H{\" o}lder-Zygmund space of order $\d$ on $M$. The following facts are well known:\\

\nd
a. ${\mc C}^\s(M)$ coincides with the Besov space $B_{\infty\infty}^\s(M)$\footnote{
See the general definition of the Besov space and the Triebel-Lizorkin space in Seeger-Sogge \cite{SS} or Triebel
\cite{Tr}. Here we follow the notations of \cite{Tr}.};\\
b. When $\s$ is a non-negative integer $k$, there exists the continuous embedding
$$H_\infty^k(M)\subsetneqq {\mc C}^k(M);$$
\nd
c. When $\s$ is a non-integer, ${\mc C}^\s(M)$ coincides with the space $C^{k,\,\d}(M)$ of H{\" o}lder
continuous functions with exponent $k+\d=\s$, where $k$ is the maximal integer $\leq\s$.\\

\nd By fact b and the asymptotic estimate \eqref{equ:ck}, we get that for a non-negative integer $k$
there holds
\begin{equation}
\label{equ:bek}
\|\chi_\l\|_{L_2\to B_{\infty\infty}^k}={\rm O}(\l^{k+(n-1)/2})
\ {\rm as}\ \l\to\infty.
\end{equation}
For a positive non-integer $\s$, there exists a unique non-negative integer $k$ and a unique $\theta\in (0,\,1)$ such that
$\s=k\theta+(k+1)(1-\theta)$.
By the complex interpolation of Besov spaces (cf Triebel \cite{Tr} or Bergh-L{\" o}fstr{\" o}m \cite{BL}),
\[[B_{\infty\infty}^k,\,B_{\infty\infty}^{k+1}]_\theta=B_{\infty\infty}^\s,\]
we obtain the estimates \eqref{equ:HZ} and \eqref{equ:Hold} in the following

\begin{thm}\label{thm:HZ}
Let $\s$ be a non-negative number, $k$ a non-negative integer and $\d\in (0,\,1)$. Let
${\mc C}^\s(M)$ and $C^{k,\,\d}(M)$ be the H{\" o}lder-Zygmund space of order $\s$ and the space of the
H{\" o}lder continuous functions with exponent $k+\d$ on $M$, respectively.
Then we have
\begin{equation}
\label{equ:HZ}
\|\chi_\l\|_{L_2\to {\mc C}^\s}={\rm O}(\l^{\s+(n-1)/2})
\ {\rm as}\ \l\to\infty.
\end{equation}
In particular, we get
\begin{equation}
\label{equ:Hold}
\|\chi_\l\|_{L_2\to C^{k,\d}}={\rm O}(\l^{k+\d+(n-1)/2})
\ {\rm as}\ \l\to\infty.
\end{equation}
Moreover,
\begin{equation}
\label{equ:HS}
\|\chi_\l\|_{L_2\to C^{\d}}=\O(\l^{\d+(n-1)/2})
\ {\rm as}\ \l\to\infty.
\end{equation}
\end{thm}

Supported by \eqref{equ:HS}, the author conjecture the "big Oh" in
the estimate \eqref{equ:HZ} can be replaced by the "Omega", that
is, there should hold that for $\s\geq 0$
\begin{equation}
\label{equ:HZS} \|\chi_\l\|_{L_2\to {\mc
C}^\s}=\O(\l^{\s+(n-1)/2}) \ {\rm as}\ \l\to\infty.
\end{equation}

We sketch the proof of \eqref{equ:HS} in the left part of this section.
Firstly, stimulated by the equivalence between \eqref{equ:c0} and \eqref{equ:spe},
we think of looking for an equivalent estimate of the square-sum of $C^\d$ norms of a bunch eigenfunctions.

\begin{lem} {\rm (\cite{Xu3})}
The validity of the estimate \eqref{equ:HS} is equivalent to the validity of the following two statements{\rm :}\\
{\rm (a)}  The estimate
\begin{equation}
\label{equ:eigenCd}
\sum_{\l_j\in (\l,\,\l+1]}\frac{\bigl(e_j(x)-e_j(y)\bigr)^2}
{\bigl({\rm dist}\,(x,\,y)\bigr)^{2\d}}= {\rm O}(\l^{(n-1)+2\d})\  {\rm as} \ \l\to\infty
\end{equation}
holds for any distinct two points $x$ and $y$ in $M$.\\
{\rm (b)} The estimate \eqref{equ:eigenCd} is sharp
in the sense that there exist a sequence of positive numbers $\{\mu_\ell\}$ with
$\lim_{\ell\to \infty}\mu_\ell=\infty$ and two sequences of points $\{x_\ell\},\,\{y_\ell\}$ on $M$
with $x_\ell\not= y_\ell$ for all $\ell$ such that
\begin{equation}
\label{equ:eigCd}
\sum_{\l_j\in (\mu_\ell,\, \, \mu_\ell+1]}\frac{\bigl(e_j(x_\ell)-e_j(y_\ell)\bigr)^2}
{\bigl({\rm dist}\,(x_\ell,\,y_\ell)\bigr)^{2\d}} \geq  C\mu_\ell^{(n-1)+2\d}.
\end{equation}
\end{lem}

Since the numerator
$$\sum_{\l_j\in (\mu_\ell,\, \, \mu_\ell+1]}\bigl(e_j(x_\ell)-e_j(y_\ell)\bigr)^2$$ on the LHS of
\eqref{equ:eigCd} at most has growth order of $\mu_\ell^{n-1}$ by \eqref{equ:spe}, if \eqref{equ:eigCd} were
true, then, to match the growth order of $\mu_\ell^{(n-1)+2\d}$ on the RHS, ${\rm dist}\,(x_\ell,\,y_\ell)$ would be
forced to converge to zero as $\ell\to\infty$ and its order should be at most of $1/\mu_\ell$. Therefore,
we can reduce the estimate \eqref{equ:eigenCd} and its sharpness \eqref{equ:eigCd} to the following estimate,
\begin{equation}
\label{equ:Omega}
\sup_{{\rm dist}(x,\,y)\sim 1/\l}\sum_{\l_j\in (\l,\,\l+1]}\frac{\bigl(e_j(x)-e_j(y)\bigr)^2}
{\bigl({\rm dist}\,(x,\,y)\bigr)^{2\d}}=\O(\l^{(n-1)+2\d})\  {\rm as} \ \l\to\infty,
\end{equation}
where ${\rm dist}\,(x,\,y)\sim 1/\l$ means that ${\rm
dist}\,(x,\,y)$ is comparable to $1/\l$. Actually we shall give an
asymptotic formula for
$$\sum_{\l_j\leq \l}\bigl(e_j(x)-e_j(y)\bigr)^2=e(x,x,\l)+e(y,y,\l)-2e(x,y,\l),$$
which implies \eqref{equ:Omega}.

\subsection{The asympotitc formula of $e(x,y,\l)|_{{\rm dist}\,(x,\,y)={\rm O}(1/\l)}$}
Consider the radial function $\Phi_n$ on ${\Bbb R}^n$,
\begin{equation*}
%\label{equ:rad}
\Phi_n(z)=\Phi_n(|z|)=(2\pi)^{-n}\int_{\{\xi\in {\Bbb R}^n:\,|\xi|\leq 1\}}
 e^{\sqrt{-1}\la z,\,\xi\ra}\,{\rm d}\xi,
\end{equation*}
where $\la\,,\,\ra$ is the inner product in ${\Bbb R}^n$. $\Phi_n(\tau)$, $\tau\geq 0$, can be rewritten as
\[\Phi_n (\tau)=(2\pi)^{-n}\,|{\Bbb S}^{n-1}|\,\int_{-1}^1 \cos\,(\tau t)(1-t^2)^{(n-1)/2}\,{\rm d}t,\]
where $|{\Bbb S}^{n-1}|$ is the area of the unit sphere ${\Bbb S}^{n-1}$ in ${\Bbb R}^n$. Note that
$\Phi(\tau)$ has a countable number of zeroes on $[0,\,\infty)$ and $\Phi_n(0)>0$. For example,
$\Phi_3(\tau)=0$ if and only if $\tan\, \tau=\tau>0$.

\begin{thm}
\label{thm:asym} {\rm (\cite{Xu3})}
Let $\tau$ be a fixed non-negative number.
As $\l\to\infty$, there holds the asymptotic formula
\begin{equation}
\label{equ:offdiag}
e(x,y,\l)|_{\l\times {\rm dist}\,(x,\,y)=\tau}=
\left\{
\begin{array}{rl}
\Phi_n(\tau) \l^n+{\rm O}(\l^{n-1}), & {\rm if}\ \Phi_n(\tau)\not=0,\\
{\rm O}(\l^{n-1}), & {\rm otherwise}.
\end{array}
\right.
\end{equation}
Moreover, as $\l\to\infty$, for any non-negative $\tau$, there holds the one term asymptotic formula
\begin{eqnarray}
\sum_{\l_j\leq \l}\Bigl(e_j(x)-e_j(y)\Bigr)^2_{\l\times {\rm dist}\,(x,\,y)=\tau}&=&
\nonumber\\
\Bigl(e(x,x,\l)+e(y,y,\l)-2e(x,y,\l)\Bigr)_{\l\times {\rm dist}\,(x,\,y)=\tau}&=&
\nonumber\\
\label{equ:asymh}
2\bigl(\Phi_n(0)-\Phi_n(\tau)\bigr)\l^n&+&{\rm O}(\l^{n-1}).
\end{eqnarray}
\end{thm}

\begin{rem} The points $x$ and $y$ in Theorem \ref{thm:asym} vary as $\l$ tends to infinity.
For two fixed distinct points $z$ and $w$ in $M$, H{\" o}rmander \cite{H1} proved that
$e(z,w,\l)={\rm O}(\l^{n-1})$ as $\l\to\infty$.
Letting $\tau=0$ in \eqref{equ:offdiag}, we get the local Weyl law \eqref{equ:asy}.

Since $\Phi_n(0)$ is always greater than $\Phi_n(\tau)$ for positive $\tau$
so that, simultaneously, the formula \eqref{equ:asymh} has a solid principal term on the RHS,
which implies \eqref{equ:Omega}.  However, the reader should be cautious at one point: Although the function
$\Phi_n(\tau)$ converges to zero as $\tau\to +\infty$ by the Riemann-Lebesgue lemma, neither
\eqref{equ:offdiag} nor \eqref{equ:asymh} make sense as $\tau=\l\,{\rm dist}\,(x,\,y)\to\infty$.

As explained in Subsection\,\ref{subsec:wk}, $e(x,y,\l)$, as ${\rm dist}\,(x,\,y)\leq {\rm inj}_M$, has local
character expressed by the Hadamard parametrix of the wave kernel. Then
we can prove Theorem \ref{thm:asym} by using the generalized
Tauberian argument (cf \cite{Xu} Lemma 2.5)  and the essential assumption that
$\l\,{\rm dist}\,(x,\,y)$ is a constant.

\end{rem}

\begin{rem} The author conjecture that there should exist similar asymptotic formulae of
$\p_x^\al \p_y^\bt e(x,y,\l)|_{\l\times {\rm dist}\,(x,\,y)=\tau}$ in a normal chart with
raduis ${\rm inj}_M$, which might be a generalization of
the asymptotic formula \eqref{equ:deri}.  If they were true, so would be the conjecture in \eqref{equ:HZS} at least
for any non-integer $\s$.
\end{rem}

\section{Fractional Sobolev $L_p$ ($2\leq p<\infty$) norms of $\chi_\l$}
With the help of the oscillatory integral theorem of Carleson-Sj{\" o}lin
\cite{CS} and Stein \cite{St0}, Sogge \cite{S2} (also see Sogge \cite{S5} for more general result)
showed that as $\l\to\infty$
\begin{equation}
\label{equ:2q}
\|{\chi}_\l\|_{L_2\to L_{2(n+1)/(n-1)}}={\rm O}\bigl(\l^{\frac{n-1}{2(n+1)}}\bigr)
\end{equation}
by using the local character of $\chi_\l$ (cf Subsection\,\ref{subsec:wk}).
Interpolating \eqref{equ:2q} with the estimate \eqref{equ:c0} and
$\|\chi_\l\|_{L_2\to L_2}=1$,
Sogge proved that
\[
||{ \chi}_\l f||_{L_2\to L_p}={\rm O}(\l^{\e(p)}),\ 2\leq p<\infty,\ {\rm as}\ \l\to\infty.
\]
where
$$\e(p)=\max\biggl(\frac{n-1}{2}-\frac{n}{p},\,\Bigl(\frac{1}{4}-\frac{1}{2p}\Bigr)(n-1)\biggr).$$
Sogge further showed by using some testing functions skillfully that
\begin{equation}
\label{equ:sog}
||{ \chi}_\l f||_{L_2\to L_p}=\O(\l^{\e(p)}),\ 2\leq p<\infty.
\end{equation}

\begin{thm}
\label{thm:Sob}
Let $s\geq 0$ and $2\leq p<\infty$.
Let $H_p^s(M)$\footnote{Cf Triebel \cite{Tr} for the definition.} be the fractional Sobolev $L_p$ space of order $s$ on $M$.
Then the following estimate hold,
\begin{equation}
\label{equ:Sob}
\|{\chi}_\l\|_{L_2\to H_p^s}=\O\,(\l^{\e(r)+s})\ {\rm as}\ \l\to\infty.
\end{equation}
\end{thm}

\nd
{\sc Proof}  Let $p'$ be the number such that $1/p'+1/p=1$. By Sogge's result \eqref{equ:sog} and the dual argument,
we have
\begin{equation*}
\|\chi_\l\|_{L_{p'}\to L_2}=\|\chi_\l\|_{L_2\to L_p}=\O(\l^{\e(p)}).
\end{equation*}
The elliptic regularity tells us that
$$\|\chi_\l \|_{L_2\to H_p^s}\sim \|(1+\D)^{s/2} \chi_\l \|_{L_2\to L_p}.$$
Since $\|(1+\D)^{s/2} \chi_\l\|_{L_{p'}\to L_2}=\|(1+\D)^{s/2} \chi_\l\|_{L_2\to L_p}$,
we only need to show
$$\|(1+\D)^{s/2} \chi_\l\|_{L_{p'}\to L_2}\sim \l^s\|\chi_\l\|_{L_{p'}\to L_2}.$$
Actually, by the self-adjointness of $(1+\D)^{s/2}$, for any $f\in L_2(M)$, there holds that as $\l\to\infty$
\begin{eqnarray*}
\|(1+\D)^{s/2} \chi_\l f\|_2^2&=&\la (1+\D)^{s/2}\chi_\l f\,,\,(1+\D)^{s/2}\chi_\l f\ra\\
&=&\la (1+\D)^s \chi_\l f\,,\,\chi_\l f\ra \sim \l^{2s} \|\chi_\l f\|_2^2,
\end{eqnarray*}
where $\la\,,\,\ra$ denotes the inner product in $L_2(M)$.
\hfill{$\Box$}

\section{${\mathcal C}^\s$ norms, $H_\infty^k$ norms and nodal domains of eigenfunctions}
Let $\vp$ be a smooth function on ${\Bbb R}$ such that
\[\vp(x)=1\ {\rm if}\ |x|\leq 1\ {\rm and}\ \vp(x)=0\ {\rm if}\ |x|\geq 3/2.\]
We put $\vp_0=\vp$, $\vp_1=\vp(x/2)-\vp(x)$ and
\[\vp_j=\vp_1(2^{1-j}\,x),\ x\in{\Bbb R},\ j\in{\Bbb N}.\]
Then $\{\vp_k\}_{k=0}^\infty$ form a dyadic resolution of unity in ${\Bbb R}$.
Seeger-Sogge (Theorem 4.1, \cite{SS}) gave characterizations of the Besov space and the Triebel-Lizorkin space
$M$ by functions $\vp_k(\D)$, $k\in {\Bbb N}_0$, of $\D$. In particular, for non-negative $\s$ there holds
\begin{equation}
\label{equ:SS}
\|f\|_{B_{\infty\infty}^\s (M)}\sim \sup_{k\in {\Bbb N}_0}\, 2^{k\s}\,\|\vp_k(\D)f\|_{L_\infty(M)},\
\forall\,f\in B_{\infty\infty}^\s (M),
\end{equation}
where "$\sim$" means the equivalence between the above two norms.

Let $e_\l$ be a real-valued non-constant eigenfunction of $\D$: $\D e_\l=\l^2 e_\l$ ($\l>0$). By the elliptic regularity,
we know that
\begin{equation}
\label{equ:SobLp}
\|e_\l\|_{H_p^s}\sim \l^s\,\|e_\l\|_{L_p}\  {\rm as}\  \l\to\infty,
\end{equation}
where $s>0$ and $1<p<\infty$.
By \eqref{equ:SS}, as $\l\to\infty$, we obtain that for $\s\geq 0$,
\begin{equation}
\label{equ:c0HZ}
\|e_\l\|_{{\mc C}^\s}\sim \sup_{k\in {\Bbb N}_0}\, 2^{k\s}\,|\vp_k(\l)|\cdot\|e_\l\|_{L_\infty}
\sim \l^\s \|e_\l\|_{L_\infty} \ {\rm as}\ \l\to\infty.
\end{equation}
By the continuous embedding $H_k^\infty\subset {\mc C}^k$, $k\in {\Bbb N}$,  it is clear
that $\l^k \|e_\l\|_{L_\infty}$ is the "big Oh" of $\|e_\l\|_{H_\infty^k}$ as $\l\to\infty$.
Furthermore, for any positive integer $k$, there holds
\begin{equation}
\label{equ:c0k}
\|e_\l\|_{H_\infty^k}\sim \l^k \|e_\l\|_{L_\infty}\ {\rm as}\ \l\to\infty.
\end{equation}
In particular, case $k=1$ of \eqref{equ:c0k} follows from a
variable coefficient version of the gradient estimate for
Poisson's equation (cf (3.15) in Gilbarg-Trudinger \cite{GT}) and
the relation \eqref{equ:c0HZ}. The justification of the general
case in \eqref{equ:c0k} only need a little bit more effort to
generalize the gradient estimate (3.15) in Gilbarg-Trudinger
\cite{GT} to derivatives of higher order. Therefore, we can
conclude that the asymptotic $L_\infty$ norm estimates for
eigenfunctions are basic in the research of eigenfunctions. We
summarize the above in the following

\begin{thm}
\label{thm:ckHZ}
Let $\s$ be a non-negative number, $k$ be a positive integer, and $e_\l$ be a real-valued non-constant
eigenfunction of eigenvalue $\l^2$, $\l>0$. Then the following equivalent relationships hold{\rm :}
\begin{equation}
\label{equ:c0kHZ}
\|e_\l\|_{{\mc C}^\s}\sim \l^\s \|e_\l\|_{L_\infty},\
\|e_\l\|_{H_\infty^k}\sim \l^k \|e_\l\|_{L_\infty}\ {\rm as}\ \l\to\infty.
\end{equation}
\end{thm}

We give an application to the geometric property of eigenfunctions as follows.

\begin{defn}
 Let $f:M\to {\Bbb R}$ be a continuous function. Then the {\it nodal set} of $f$ is the set $f^{-1}(0)$,
the {\it nodal domain} of $f$ is a connected component on $M\backslash f^{-1}(0)$, and the {\it concentrating set} of
$f$ is the set of points  where $f$ attains its maximum or minimum.
\end{defn}

Let $\O$ be an open, bounded domain in the Riemannian manifold $(M^n,\,g)$.
If $x\in \O$, then we denote $d(x)$ is the minimum distance from $x$ to the boundary
$\p \O$. We denote by $R(\O)=\max_{x\in\O}\, d(x)$ the inner radius of $\O$, that is, the radius of the largest inscribed
ball.

\begin{thm}
\label{thm:nodal}
Let $e_\l$ be a real-valued non-constant eigenfunction of eigenvalue $\l^2$, $\l>0$.
Then the distance between the nodal set and the concentrating set of $e_\l$ is comparable to $1/\l$ as $\l\to\infty$.
Moreover, there exists at least two nodal domains of $e_\l$ whose inner radii are comparable to $1/\l$ as $\l\to\infty$.
\end{thm}

\begin{rem}
\label{rem:nodal} Let $\O$ be any nodal domain of $e_\l$. It is
known (cf Lemma 10 in Savo\footnote{Savo only proved the dimension
two case. But his proof is also valid to general dimension with a
little modification.} \cite{Savo}) that $R(\O)\leq C/\l$, where
the constant $C$ is depending on the lower bound of Ricci tensor
of $M$ and the diameter of $M$. As a consequence, for any point
$x\in M$, the minimum distance of $x$ to the nodal set of $e_\l$
is $\leq C/\l$. However, it seems difficult to give a upper bound
for the inner radius of an arbitrary nodal domain of $e_\l$.
\end{rem}

\nd{\sc Proof of Theorem \ref{thm:nodal}}\quad Let $z$ and $w$ belong to in the nodal set and the concentrating set of
$e_\l$, respectively. By Remark \ref{rem:nodal}, to prove the first statement, we need only to prove that the distance between $z$ and $w$ is
$\geq C/\l$ for a constant independent of $\l$. By Theorem 2 in Nadirashivili \cite{Na}, $\max_M e_\l/|\min_M e_\l|$ is
comparable to 1 as $\l\to\infty$. We may assume that $e_\l$ takes  value $\|e_\l\|_{L_\infty}$ at $w$.
Let $\d$ be a number in $(0,\,1)$. Then
\[\frac{\|e_\l\|_{L_\infty}}{\bigl({\rm dist}\,(z,\,w)\bigr)^\d}=
\frac{|e_\l(w)-e_\l(z)|}{\bigl({\rm dist}\,(z,\,w)\bigr)^\d}\leq
\|e_\l\|_{C^\d}.\]
Hence, by Theorem \ref{thm:ckHZ} we get ${\rm dist}\,(z,\,w)\geq C/\l$.

Let $w_1$ and $w_2$ be the points where $e_\l$ attains its maximum and minimum, respectively.
Theorem 2 in Nadirashivili \cite{Na} tells us that both $e_\l(w_1)$ and $e_\l(w_2)$ are comparable to
the $L_\infty$ norm of $e_\l$. Let $\O_1$ and $\O_2$ be two nodal domains which contains $w_1$ and $w_2$, respectively.
By the previous result, ${\rm dist}\,(w_k,\,\p \O_k)\geq C/\l$ for $k=1,2$. Hence,
the two nodal domains $\O_1$ and $\O_2$ have inner radii comparable to $1/\l$.\hfill{$\Box$}

\section{Some spherical harmonics}
Let $k$ be a non-negative integer and $2\leq p<\infty$, $\s\geq 0$.
Let $e_\l$ be a non-constant eigenfunction: $\D e_\l=\l^2 e_\l$, $\l>0$.
The estimates \eqref{equ:ck}, \eqref{equ:HZ} and \eqref{equ:Sob} of $\chi_\l$ give the corresponding
asymptotic upper bounds for $e_\l$ as $\l\to\infty$:
\begin{eqnarray}
\label{equ:1}
\|e_\l\|_{H_\infty^k}/\|e_\l\|_2&=&{\rm O}(\l^{k+(n-1)/2}),\\
\label{equ:2}
\|e_\l\|_{{\mc C}^\s}/\|e_\l\|_2&=&{\rm O}(\l^{\s+(n-1)/2}),\\
\label{equ:3}
\|e_\l\|_{H_p^s}/\|e_\l\|_2&=& {\rm O}(\l^{s+\e(p)}).
\end{eqnarray}
Let ${\Bbb S}^n$ be the unit sphere of dimension $n$ in ${\Bbb R}^{n+1}$.
Let $Z_m$ be the zonal harmonic function of degree $m$ with respect to the north pole and $Q_m$
be the spherical harmonics defined by $Q_m(x)=(x_2+\sqrt{-1} x_1)^m$. Remember that both $Z_m$ and $Q_m$ are eigenfunctions
of eigenvalue $m(m+n-1)$ on ${\Bbb S}^n$.  Then Sogge \cite{S0} showed
that there exists a positive constant $C$ independent of $m$ such that
\begin{eqnarray*}
\|Z_m\|_r/\|Z_m\|_2\geq C m^{\e(r)}, & & 2(n+1)/(n-1)\leq r\leq \infty,\\
\|Q_m\|_r/\|Q_m\|_2\geq C m^{\e(r)}, & & 2\leq r\leq 2(n+1)/(n-1).
\end{eqnarray*}
That is, the asymptotic upper bounds for $L_r$ norm of eigenfunctions are attained by
$Z_m$ for $2(n+1)/(n-1)\leq r\leq \infty$ and $Q_m$ for $2\leq r\leq 2(n+1)/(n-1)$. By \eqref{equ:SobLp}, it is also the
case for the upper bounds in \eqref{equ:3}. By Theorem \ref{thm:ckHZ}, the zonal harmonic functions $Z_m$ attain the upper
bounds \eqref{equ:1} and \eqref{equ:2}.  \\

\nd {\bf Acknowledgements} Part of this research has been done
through the winter of 2004 and the spring of 2005, when I visited
Tokyo Institute of Technology and Johns Hopkins University,
respectively. I thank the departments of mathematics of these two
universities for their warm hospitality. I also express my great
gratitude to Professors C. D. Sogge, S. Zelditch and X.-J. Xu for
their interest in my work and their valuable conversations.
Special thanks goes to Professors T. Ochiai and H. Arai for their
constant encouragement. At last, I appreciate the careful reading
of the manuscript by the referee and her/his valuable comments.
This research was supported in part by the JSPS Postdoctoral
Fellowship for Foreign Researchers, the FRG Postdoctoral
Fellowship and the Program of Visiting Scholars at Nankai
Institute of Mathematics.

%%%==========================================================


\begin{thebibliography}{99}

\bibitem{A}{V. G. Avakumovi${\rm {\check c}}$,}
{\it {\" U}ber die eigenfunktionen auf geschlossenen Riemannschen mannigfaltigkeiten},
Math. Z. {\bf 65} (1956), 327-344.

\bibitem{BB1}{R. Balian and C. Bloch},
{\it Distribution of eigenfrequencies for the wave equation in a finite domain I:
Three-dimensional problem with smooth boundary surface},  Ann. Phys. {\bf 60}
(1970), 401-447.

\bibitem{BB2}{R. Balian and C. Bloch},
{\it Distribution of eigenfrequencies for the wave equation in a finite domain II:
Electromagnetic field, Riemannian spaces}, Ann. Phys. {\bf 64}
(1971), 271-307.


\bibitem{BB3}{R. Balian and C. Bloch},
{\it Distribution of eigenfrequencies for the wave equation in a finite domain III:
Eigenfrequency density oscillations}, Ann. Phys. {\bf 69}
(1972), 76-160.

\bibitem{BL} J. Bergh and J. L{\" o}fstr{\" o}m,
{\it Interpolation Spaces: An Introduction},
Springer-Verlag, Berlin, Heidelberg, 1976.

\bibitem{CS}L. Carleson and P. Sj{\" o}lin,
{\it Oscillatory Integrals and a Multiplier Problem for the Disc},
Studia. Math. {\bf 44} (1972), 287-299.

\bibitem{GT} D. Gilbarg and N. S. Trudinger,
{\it Elliptic Partial Differential Equations of Second Order},
Springer-Verlag, Berlin, Heidelberg, New York, 2001.

\bibitem{G}D. Grieser,
{\it Uniform Bounds for Eigenfunctions of the Laplacian on Manifolds
with Boundary}, Comm. Partial Differential Equations. {\bf 27} (2002), 1283-1299.



\bibitem{H1}L. H{\" o}rmander,
{\it The spectral function of an elliptic operator},
Acta Math. {\bf 88} (1968), 341-370.


\bibitem{H2}L. H{\" o}rmander,
{\it The Analysis of Linear Partial Differential Equations III},
Corrected second printing, Springer-Verlag, Tokyo 1994;

\bibitem{Lan} L. D. Landau and E. M. Lifshitz,
{\it Theory of Elasticity (Volume 7 of Course of Theoretical Physics)},
Pergamon Press London, 1959.

\bibitem{Lev} B. M. Levitan,
{\it On the asymptotic behaviour of the spectral function of a self-adjoint differential equation of
second order}, Isv. Akad. Nauk SSSR Ser. Mat. {\bf 16} (1952), 325-352.

\bibitem{Na} N. S. Nadirashvili,
{\it Metric properties of eigenfunctions of the Laplace operator on manifolds},
Ann. Inst. Fourier, Grenoble, {\rm 41} (1991), 259-265.

\bibitem{P1}{$\stackrel{\circ}{\rm A}$. Pleijel,}
{\it A study of certain Green's function with applications in the theory of vibrating membranes},
Ark. Math. {\bf 2} (1954),  553-569.

\bibitem{P2}{$\stackrel{\circ}{\rm A}$. Pleijel,}
{\it On Green's function and the eigenvalue-distribution of the three dimensional membrane equation},
Tolfte Scand. Mat. Kongr. Lund, 1953, pp. 222-240. Lunds Universitets Mat. Inst. Lund (1954).


\bibitem{Sa} Yu. Safarov and D. Vassiliev,
{\it The Asymptotic Distribution of Eigenvalues of Partial Differential Operators,}
Translations of Mathematical Monographs, {\bf 155},
American Mathematical Society, Rhode Island, 1997.

\bibitem{Sar} P. Sarnak,
{\it Arithmetic quantum chaos}, Israel Math. Conf. Proc. {\bf 8},
Bar-Ilan Univ., Ramat Gan, 1995, 183-236.

\bibitem{Savo} A. Savo,
{\it Lower bounds for the nodal length of eigenfunctions of the Laplacian},
Ann. Glob. Anal. Geo. {\bf 19} (2001), 133-151.

\bibitem{SS} A. Seeger and C. D. Sogge,
{\it On the boundedness of functions of (pseuo)differential operators on compact manifolds},
Duke Math. J. {\bf 59} (1989), 709-736.

\bibitem{S0}Sogge, C. D.
{\it Oscillatory integrals and spherical harmonics},
Duke Math. J.  {\bf 53} (1986), 43-65.


\bibitem{S1}C. D. Sogge,
{\it Concerning the $L^p$ norm of spectral clusters for second order elliptic
operators on compact manifolds},
J. Funct. Anal.  {\bf 77} (1988),  123-134.

\bibitem{S2}C. D. Sogge,
{\it Remarks on $L^2$ restriction theorems for Riemannian manifolds}, in:
E. Berkson, T. Peck and J. Jerry Uhl, Jr. (eds),
{\it Analysis at Urbana},
London Math. Soc. Lecture Note Ser. {\bf 137},
Cambridge Univ. Press, Cambridge, 1989; Vol. 1, 416-422.


\bibitem{S3} C. D. Sogge,
{\it On the convergence of Riesz means of compact manifolds},
Ann. Math.  {\bf 126} (1987), 439-447.

\bibitem{S5} C. D. Sogge,
{\it Fourier Integrals in Classical Analysis}, Cambridge University Press, 1993.

\bibitem{S4}C. D. Sogge,
{\it Eigenfunction and Bochner-Riesz estimates on manifolds with boundary},
Mathematical Research Letter,  {\bf 9} (2002), 205-216.

\bibitem{SZ} C. D. Sogge and S. Zelditch,
{\it Riemannian manifolds with maximal eigenfunction growth},
Duke Math. J. {\bf 114} (2002), 387-437.

\bibitem{St0}E. M. Stein,
{\it Oscillatory integrals in Fourier analysis}, in: E. M. Stein (eds),
Beijing Lectures in Harmonic Analysis,
Princeton Univ. Press: Princeton New Jersey, 1986; 307-356.


\bibitem{Tr} H. Triebel,
{\it Theory of Function Spaces},
Birkh{\" a}user Verlag, Basel, Boston, Stuttgart, 1983.

%\bibitem{Sz} Szeg{\" o}, G.
%{\it Orthogonal Polynomials.}
%American Mathematical Society Colloquium Publications, {\bf 23},
%American Mathematical Society, New York, 1939.

\bibitem{TZ} J. A. Toth and S. Zelditch,
{\it Riemannian manifolds with uniformly bounded eigenfunctions},
Duke Math. J. {\bf 111} (2002), 97-132.

\bibitem{Xu} B. Xu,  {\it Derivative of the spectral function and Sobolev norms of eigenfunctions on a
closed Riemannian manifold}, Ann. Glob. Anal. Geo. {\bf 26} (2004), 231-252.

\bibitem{Xu2} B. Xu,
{\it Derivative of spectral function and Sobolev norms of eigenfunctions on a closed Riemannian manifold},
Harmonic Analysis and Nonlinear PDEs, Surikaisekikenkyusho Kokyuroku, {\bf 1389} (2004), 60-77.

\bibitem{Xu3} B. Xu,
{\it Concerning the H{\" o}lder-Zygmund Norm of Spectral Clusters on a Closed Riemannian Manifold
}, in preparation.

\bibitem{X} X.-J. Xu,
{\it Eigenfunctions function estimates on compact Riemannian manifolds
with boundary and H{\" o}rmander multiplier theorem}, Thesis, Johns Hopkins University.

\bibitem{Z1} S. Zelditch,
{\it Uniform distribution of eigenfunctions on compact hyperbolic surfaces},
Duke Math. J. (1987), 919-941.



\end{thebibliography}
\end{document}